\documentclass{amsart}

\usepackage{amsmath,amssymb}
\newtheorem{theorem}{Theorem}[section]
\newtheorem{lemma}[theorem]{Lemma}
\newtheorem{proposition}[theorem]{Proposition}

\newtheorem{example}[theorem]{Example}

\begin{document}

\title{Modules with RD-composition
  series over a commutative ring}  

\author{Fran\c cois  Couchot}

\begin{abstract} If $R$ is a commutative ring, then we prove that every
finitely generated $R$-module has a pure-composition series with
indecomposable cyclic factors and any two such series are isomorphic if and
only if $R$ is a B\' ezout ring and a CF-ring. When $R$ is a such ring, the
length of a pure-composition series of a finitely generated $R$-module $M$
is compared with its Goldie dimension and we prove that these numbers are
equal if and only if $M$ is a direct sum of cyclic modules. We also give an
example of an artinian module over a noetherian domain, which has an RD-
composition series with uniserial factors. Finally we prove that every
pure-injective $R$-module is RD-injective if and only if $R$ is an
arithmetic ring.
\end{abstract}
\maketitle

\bigskip
In this paper, for a commutative ring $R,$ we study the following
properties:
\begin{enumerate}
\item Every finitely generated $R$-module $M$ has a finite chain of
  RD-submodules with cyclic factors.
\item Every finitely generated $R$-module $M$ has a finite chain of
  RD-submodules with indecomposable cyclic factors.
\item $R$ satisfies (2) and any two chains of RD-submodules of $M,$
  with indecomposable cyclic factors, are isomorphic.
\end{enumerate}

In \cite{SaZa}, L.Salce and P.Zanardo proved that every valuation ring
satisfies (3), and in \cite{Nau} C.Naud\' e showed that each h-local
and B\' ezout domain also satisfies (3).

In section~\ref{S:pre}, we give definitions and some preliminary results. It is
proved that every ring that satisfies (1), is a Kaplansky ring(or an
elementary divisor ring), but we don't know if the converse holds.

In section~\ref{S:ind}, we state that a ring $R$ satisfies (3) if and only if
$R$ is a B\' ezout ring and a CF-ring(\cite{ShWi}). We show that $R$
satisfies (2), when $R$ is a semilocal arithmetic ring or an
h-semilocal B\' ezout domain, and for every finitely generated
$R$-module $M,$ we prove that any chain of RD-submodules of $M$ has
the same length which is equal to the number of terms of a finite
sequence of prime ideals, associated to $M.$ However, we give two
examples of semilocal arithmetic rings that don't satisfy (3).

In section~\ref{S:gol}, when $R$ is a ring that verifies (3), the length $l(M)$
of each chain of RD-submodules, with indecomposable cyclic factors, is
compared with $g(M),$ the Goldie dimension of $M.$ We show that $g(M)\leq
l(M)$ and that $M$ contains a direct sum of $g(M)$ nonzero
indecomposable cyclic submodules which is an essential RD-submodule of
$M.$ Hence it follows that $M$ is a direct sum of cyclic submodules if
and only if $g(M)=l(M).$ These results were proved in \cite{SaZa}, when $R$ is
a valuation ring and in \cite{Nau}, when $R$ is an h-local B\' ezout
domain.

Finally, in section~\ref{S:inj}, we give an example of an artinian
module, over a noetherian domain, that has a finite chain of
RD-submodules, with uniserial factors. We also give explicit examples
of pure-injective modules that fail to be RD-injective, over
noetherian domains. We prove that every pure-injective $R$-module is
RD-injective if and only if $R$ is an arithmetic ring. When $R$ is a
domain, this result was proved in \cite{NaPr}, by C.G Naud\' e, G. Naud\' e and
L.M. Pretorius. 

\section{Preliminaries}
\label{S:pre}
All rings in this paper are commutative with unity, and all modules are
unital. We recall that a module, over a ring $R,$ is
\textit{uniserial} if its set of submodules is totally ordered by
inclusion. A ring $R$ is a \textit{valuation ring} if $R$ is a
uniserial module, and $R$ is \textit{arithmetic} if $R_P$ is a
valuation ring for every maximal ideal $P.$ We say that $R$ is a
\textit{B\' ezout ring} if every finitely generated ideal is principal
and $R$ is a \textit{Kaplansky ring (or elementary divisor ring)} if
every finitely presented $R$-module is a finite direct sum of cyclic
presented modules (see \cite{Kap} and \cite{LLS}). 

An exact sequence of $R$-modules : $0\rightarrow F\rightarrow
E\rightarrow G\rightarrow 0$ is \textit{pure-exact} if it remains exact when
tensoring it with any $R$-module. In this case we say that $F$ is a
\textit{pure submodule} of $E.$
When $rE\cap F = rF$ for every $r\in R,$ we say that $F$ is an
\textit{RD-submodule} of $E$(relatively divi\-sible) and that the sequence
is \textit{RD-exact.}
 Then every pure submodule is an RD-submodule, and the equivalence
 holds if and only if $R$ is an arithmetic ring (\cite{War1} and \cite{War2}). We say
 that an $R$-module $M$ has a \textit{pure-composition series}
 (respectively an \textit{RD-composition series})
if there exists a finite chain $\{0\}=M_0\subset M_1\subset\dots\subset M_n = M$ of
pure (respectively RD-) submodules. If we denote $A_i =
$ann$(M_i/M_{i-1}),$ then $(A_i)_{1\leq
i\leq n}$ is called the \textit{annihilator sequence} of the pure-composition
series. We say that the annihilator sequence is \textit{increasing} if $A_1\subseteq
A_2\subseteq\dots\subseteq A_n,$ and it is \textit{totally ordered} if
$A_i$ and $A_j$ are comparable $\forall i, \forall j, 1\leq i,j\leq n.$

An $R$-module $F$ is \textit{pure-projective} (respectively
\textit{RD-projective}) if for every pure (respectively RD-) exact
sequence: 
 $0\rightarrow N\rightarrow M\rightarrow L\rightarrow 0,$ of
$R$-modules, the following sequence:  
$0\rightarrow\hbox{Hom}_R(F,N)\rightarrow
\hbox{Hom}_R(F,M)\rightarrow\hbox{Hom}_R(F,L)\rightarrow 0,$ is
exact.

An $R$-module $F$ is \textit{pure-injective} (respectively
\textit{RD-injective}) if for every pure (respectively RD-) exact
sequence: 
$0\rightarrow N\rightarrow M\rightarrow L\rightarrow 0$ of
$R$-modules, the following sequence: 
 \(0\rightarrow\hbox{Hom}_R(L,F)
\rightarrow\hbox{Hom}_R(M,F)\rightarrow\hbox{Hom}_R(N,F)\rightarrow
0,\) is exact.  
Then we have the following proposition.
\begin{proposition} 
\label{P:Kap}
Let $R$ be a ring. The following conditions
are equivalent:

\begin{enumerate}
\item Every finitely generated module has a pure-composition series
with cyclic factors.

\item Every finitely generated  module has an RD-composition
series with cyclic factors.
\end{enumerate}
Moreover, if $R$
satisfies one of these conditions, $R$ is a Kaplansky ring.
\end{proposition}
\textbf{Proof.} If $R$ satifies (1), by the same proof as in \cite[Theorem
2.3 and Corollary 2.4]{Nau}, we show that $R$ is a Kaplansky ring.

Since every pure submodule is an RD-submodule, it follows that
$(1)\Rightarrow (2).$

If $R$ satisfies (2) and if we prove that $R$ is an arithmetic ring,
by \cite[Theorem 3]{War1} $R$ satisfies~(1).

For every maximal ideal $P,$ we must prove that $R_P$  is a valuation ring.
Clearly $R_P$ also satisfies~(2). We may assume that $R$ is local and that
$P$ is its maximal ideal.
If $R$ is not a valuation ring, then, as in the proof of \cite[Theorem
2]{War2}, we may suppose that there exist $a$ and $b$ in $P$ such that $Ra\cap Rb =\{0\}$ and $\hbox{ann}(a) = \hbox{ann}(b) = P.$ Moreover, there exists an
indecomposable $R$-module $M,$ with two generators $x$ and $y,$ satisfying
the relation $bx-ay = 0.$ Since $R$ satisfies (2), there exists an
RD-composition series: $\{0\}=M_0\subset M_1\subset\dots\subset M_n = M,$
such that $M_i/M_{i-1}= R(z_i+M_{i-1})$ for each
$i,\ 1\leq i\leq n.$ For every $i,$ there exist $s_i$ and $t_i$ in
$R$ such that $z_i = s_ix+t_iy.$ Clearly $M$ is ge\-nerated by
$\{z_i\mid 1\leq i\leq n\}.$ Since $\dim_{R/P} (M/PM)= 2,$ there exist at least 2 indexes $j$ and $k,\ 1\leq j<k\leq
n$ such that $(s_j,t_j)\notin P\times P$ and $(s_k,t_k)\notin P\times
P.$ Let $j$ be the smallest index such that $(s_j,t_j)\notin P\times P.$
Then $\dim_{R/P} (M_j/PM_j)= 1,$ whence $M_j =
Rz_{j}.$ We may assume that there exists $z\in M,\ z = x+ty,$ where
$t\in R,$ such that $Rz$ is an RD-submodule of $M.$ Then $M/Rz = R(y+Rz)$ and if $I = \hbox{ann}(M/Rz),\ I = \{r\mid ry\in Rz\}.$
Let $r\in I.$ There exists $c\in R$ such that $ry = c(x+ty).$ Since
$cx+(tc-r)y=0,$ it follows that $\exists u\in R$ such that $c=ub$ and
$r-tc=ua,$ whence $r = u(a+tb).$ Therefore $I\subseteq R(a+tb)\subset
Ra\oplus Rb$ which is semi-simple. We deduce that $I = Rd$ where $d=0$ or
$d=a+tb.$ We get that $M/Rz \simeq R/Rd$ is RD-projective. Consequently
$Rz$ is a direct summand of $M,$ whence a contradiction. \qed

\bigskip
We don't know if the converse of the second assertion of this proposition holds. The following
propositions state that some classes of arithmetic rings satisfie the
equivalent conditions of the Proposition~\ref{P:Kap}, and therefore these
rings are Kaplansky.

We begin with the classe of semi-local arithmetic rings, and we
already know
that these rings are Kaplansky by \cite[Corollary 2.3]{LLS}.

If $M$ is a finitely generated $R$-module, we denote $\mu(M)$ the minimal
number of generators of $M.$
\begin{proposition} 
\label{P:semlo}
Let $R$ be a semilocal arithmetic
ring. Then every finitely generated $R$-module $M$ has a
pure-composition series with cyclic factors and an increasing annihilator sequence.
\end{proposition}
\textbf{Proof.}
Let $P_1,\dots,P_n$ be the maximal ideals of $R,$ and $J$ its Jacobson
radical. As in the proof of \cite[Theorem 15]{FuSa1}, we state that,
for
 every $k,\
1\leq k\leq n,$ there exists $x_k\in M$ such that, for each 
$y\in
x_k+P_kM_{P_k},\ R_{P_k}y$ is a pure submodule of
$M_{P_k},$ and
 $\hbox{ann}_{R_{P_k}}(R_{P_k}y) =
\hbox{ann}_{R_{P_k}}(M_{P_k}).$ Moreover $\mu_{R_{P_k}}(M_{P_k}/R_{P_k}y) =
\mu_{R_{P_k}}(M_{P_k})-1$ if $M_{P_k}\not=\{0\}.$ When
$M_{P_k}=\{0\}$ we put $x_k=0.$ By the chinese remainder Theorem,
we get $M/JM\simeq \prod_{k=1}^n(M/P_kM)
\simeq\prod_{k=1}^n (M_{P_k}/P_kM_{P_k}).$ Consequently there exists
$x\in M$ such that $x\equiv x_k \mod P_k M_{P_k},\ \forall k,\ 1\leq k\leq
n.$ We deduce that $Rx$ is a pure submodule of $M,$ that $\hbox{ann}(Rx) =
\hbox{ann } M$ and that $\mu(M/Rx)
=\mu(M)-1.$ We complete the proof by induction on $m = \mu(M).$ \qed

\bigskip
A domain $R$ is said to be {\it h-semilocal} if $R/I$ is semilocal for every
nonzero ideal $I,$ and $R$ is said to be {\it h-local} if, in addition, $R/P$
is local for every nonzero prime ideal $P,$ \cite{Mat}.

\begin{proposition} 
\label{P:mipr}
Let $R$ be a B\' ezout ring with a unique
minimal prime ideal $Q.$ We assume that $Q$ is a uniserial module and
 $R/Q$ an h-semilocal domain. Then every finitely generated
$R$-module has a pure-composition series with cyclic factors and 
an increasing annihilator sequence.
\end{proposition}
\textbf{Proof.}
We may assume that $R$ is not semilocal. If $Q\not=\{0\},$ then from
\cite[Lemma 17 and 18]{WiWi} there exists only one maximal ideal $P$ such that
$QR_P\not=\{0\},$ every ideal contained in $Q$ is comparable to every ideal
of $R$, $Q^2 = \{0\}$ and $Q$ is a torsion divisible $R/Q$-module. Let $M$ be a finitely generated module. Since
$R/Q$ is a B\'ezout domain, $M/QM\simeq F\oplus
T,$ where $T$ is a torsion $R/Q$-module and $F$ a free $R/Q$-module of
finite rank. Let $m$ be this rank. We prove this proposition by induction on $m.$

If $m=0,$ then for every $x\in M,$ there exists $s\in R\setminus Q$ such
that
$sx\in QM.$ Since $Q$ is uniserial, there exists $t\in Q$ and $y\in M$ such
that
$sx = ty.$ Since $Q$ is a torsion $R/Q$-module there exists $r\notin Q$ such
that $rt=0.$ Since $M$ is finitely generated, we get that
$Q\subset\hbox{ann}(M),$ and that $R/\hbox{ann}(M)$ is semilocal.
From the previous proposition we deduce the result.

Now suppose $m\geq 1.$ Let $\{x_1+QM,\dots,x_m+QM\}$ be a basis of
$F,$ and $\{y_1+QM,\dots,y_p+QM\}$ a spanning set of $T.$ Since $Q$ is the nilradical of $R,$ then, by Nakayama Lemma
$\{x_1,\dots,x_m\}\cup\{y_1,\dots,y_p\}$ generates $M.$ Let \(H = \{j \mid
1\leq j\leq m, \hbox{ann}(x_j) = \hbox{ann}(M)\}.\) Since $P$ is the only
maximal ideal such that $QR_P\not=\{0\},$ then \(H = \{j \mid
1\leq j\leq m, \hbox{ann}_{R_P}(R_Px_j) = \hbox{ann}_{R_P}(M_P)\}\)
also. From
\cite[Theorem 15]{FuSa1}, it follows that there exists $j\in H$ such that $R_Px_j$
is a pure-submodule of $M_P.$ Let $N$ be a maximal ideal, $N\not= P.$ Then
$M_N\simeq (M/QM)_N = F_N\oplus T_N.$
Consequently $R_Nx_j$ is a direct summand of $M_N.$ We deduce that $Rx_j$
is a pure-submodule of $M$ and if $M' = M/Rx_j,$ then
$M'/QM'\simeq F'\oplus T$ where $F'$ is a free
$R/Q$-module of rank $(m-1).$ From the induction hypothesis it follows that
$M$ has a pure-composition series with cyclic factors and an increasing
annihilator sequence. \qed

\begin{proposition} 
\label{P:pugen}
Let $R$ be an arithmetical ring and $M$ be a
finitely generated $R$-module. If $M$ has a pure-composition series with
cyclic factors and an increasing annihilator sequence, then $\mu(M)$ is the
length of this pure-composition series.
\end{proposition}
\textbf{Proof.}
Let $\{0\}=M_0\subset M_1\subset M_2\dots\subset M_{n-1}\subset M_n = M$ be a
pure-composition series such that $A_i =
\hbox{ann}(M_i/M_{i-1})\subseteq A_{i+1}$ for
every $i,\ 1\leq i\leq n-1,$ and such that $M_i/M_{i-1}
= R(x_i+M_{i-1})$ for every $i,\ 1\leq i\leq n.$ Then, clearly $\{x_i\mid
1\leq i\leq n\}$ generates $M.$ Therefore $n\geq\mu(M).$ Since
$M/M_{n-1}\not=\{0\},$ there exists a maximal ideal $P$
such that $(M/M_{n-1})_P\not=\{0\}.$
Consequently
$(M_i/M_{i-1})_P\not=\{0\}$ for every $i,\
1\leq i\leq n.$ Then $M_P$ has a pure-composition series of length $n.$ From
\cite[Lemma 1.4]{SaZa}  we deduce that $n=\mu_{R_P}(M_P).$ But
$\mu_{R_P}(M_P)\leq\mu(M)$ and consequently $n=\mu(M).$ \qed

\bigskip
Now, as in \cite{SaZa} or \cite{FuSa2}, we establish the isomorphy of any two
pure-composition series, with cyclic factors and totally ordered
annihilator sequences, i.e. the existence of a bijection between the two sets of
cyclic factors, such that corresponding factors are isomorphic. We use
similar lemmas, with the same proofs, and we get the following theorem:
\begin{theorem} 
\label{T:isord}
Let $R$ be an arithmetic ring. Then any two
pure-composition series of a finitely generated $R$-module, with cyclic
factors and totally ordered annihilator sequences, are isomorphic.
\end{theorem}
\textbf{Proof.}
See the proof of \cite[Theorem 1.6, p. 177]{FuSa2}. \qed

\section{Pure-composition series with indecomposable cyclic factors}
\label{S:ind}
We follow T.S. Shores and R. Wiegand \cite{ShWi}, by defining a {\it canonical form} for
an $R$-module $M$ to be a decomposition $M\simeq R/I_1\oplus
R/I_2\oplus\dots\oplus R/I_n,$ where $I_1\subseteq
I_2\subseteq\dots\subseteq I_n\not= R,$ and by calling a ring $R$ a
{\it CF-ring} if every direct sum of finitely many cyclic modules has a
canonical form. Now, following P. V\'amos \cite{Vam}, we say that $R$
is a {\it torch ring} if the following conditions are satisfied~:
\begin{enumerate}
\item $R$ is a nonlocal arithmetic ring.
\item $R$ has a unique minimal prime ideal $Q$ which is a nonzero uniserial
$R$-module.
\item $R/Q$ is an h-local domain.
\end{enumerate}

We will say that $R$ is a {\it semi-torch ring} if $R$ satisfies the
conditions (1), (2) and (3'): $R/Q$ is an h-semilocal domain.

By \cite[Theorem 3.12]{ShWi}, a CF-ring is a finite product of
indecomposable arithmetic rings. If $R$ is an indecomposable CF-ring,
then $R$ is a valuation ring, or an h-local domain, or a torch ring.

We will say that $R$ is a {\it semi-CF-ring} if $R$ is a finite product of 
indecomposable arithmetic rings, where each indecomposable factor ring
is semilocal, or an h-semilocal domain, or a semi-torch ring.

In the sequel, we will call a ring $R$ to be a {\it PCS-ring} if every
finitely generated module has a pure-composition series with indecomposable
cyclic factors and if any two such pure-composition series are isomorphic.

We will state the following theorem which is one of the main results
of this paper. 
\begin{theorem} 
\label{T:main}
Let $R$ be a ring. Then the following assertions are
equivalent~:
\begin{enumerate}
\item $R$ is a B\'ezout ring and a CF-ring.
\item $R$ is a PCS-ring.
\end{enumerate}
\end{theorem}

Some preliminary results are needed to prove this theorem. First we
will state some results on finitely generated modules over a B\' ezout
semi-CF-ring $R$. Recall that a module $M$ has {\it Goldie dimension} $n$
if it has an independant set of uniform submodules $M_1,\dots,M_n$ such that $M_1\oplus\dots\oplus M_n$
is essential in $M.$ We will denote $g(M)$ the Goldie dimension of $M.$
Then $g(M) = n$ if and only if the injective hull $E(M)$ of $M$ is a direct
sum of $n$ indecomposable injective modules.

\begin{proposition} 
\label{P:pugo}
Let $R$ be a B\'ezout semi-CF-ring. Then every finitely generated $R$-module has a pure-composition series with
indecomposable cyclic factors and a finite Goldie dimension.
\end{proposition}
\textbf{Proof.}
Let $M$ be a finitely generated $R$-module. We put $R=\prod_{j=1}^m
R_j$, where $R_j$ is indecomposable. Then $M\simeq \oplus_{j=1}^m
M_j$, where $M_j=R_j\otimes_R M.$ From Propositions~\ref{P:semlo}
and~\ref{P:mipr}, $M_j$ has a  pure-composition series: \[0\subset
M_{1,j}\subset\dots\subset M_{n_j,j} = M_j,\] with cyclic factors.
Then the chain (c):
\[0\subset M_{1,1}\subset\dots\subset M_{n_1,1} = M_1\subset
M_1\oplus M_{1,2}\subset\dots\subset M_1\oplus M_2\subset\dots\]
\[\dots\subset
M_1\oplus\dots\oplus M_{m-1}\oplus M_{1,m}\subset\dots\subset M_1\oplus\dots\oplus M_m = M\] is a
pure-composition series of $M,$ with cyclic factors. Let $N$ be a
factor of this pure-composition series. Then $N$ is a module over an
indecomposable factor ring $R'$ of $R.$ When $R'$ is a domain and
$\hbox{ann}_{R'}(N)=\{0\}$, or when $R'$ is a semitorch ring and
$\hbox{ann}_{R'}(N)\subseteq Q,$ where $Q$ is the minimal prime of $R',$ $N$
is indecomposable. In the other cases, $R/\hbox{ann}(N)$ is semilocal. Since
each semilocal arithmetic ring has only finitely many minimal prime
ideals which are pairwise comaximal, it follows that $N$ is a finite direct sum
of indecomposable cyclic modules. Consequently, we deduce from (c) a
pure-composition series of $M,$ with indecomposable cyclic factors.

It is sufficient to prove that $M_j$ has a finite Goldie dimension for
every $j,$ $1\leq j\leq m.$ We may assume that $R$ is indecomposable. First, we suppose that $M$ is cyclic. When $R$ is a domain and
$\hbox{ann}(M)=\{0\}$, or when $R$ is a semitorch ring and
$\hbox{ann}(M)\subseteq Q,$ where $Q$ is the minimal prime of $R,$ we have
$g(M)=1.$ In the other cases
$R/\hbox{ann}(M)$ is a semilocal ring. We may assume that $R$ is
semilocal. Let $(S_i)_{i\in I}$ be a family of independant submodules
of $M$ such that $\oplus_{i\in I}S_i$ is essential in $M.$ For
every maximal ideal $P$, $M_P$ is a uniserial $R_P$-module and
consequently there is at most one index $i$ in $I$ such that
$(S_i)_P\not=\{0\}.$ It follows that $I$ is a finite set. Hence
$M$ has a finite Goldie dimension. In the general case, let $N$ be a submodule of $M$ such that
$\mu (N)=\mu (M)-1$ and $M/N$ is cyclic.The inclusion map $N\rightarrow
E(N)$ can be extended to $w:M\rightarrow E(N).$  Let
$f:M\rightarrow E(N)\oplus E(M/N)$ defined by
$f(x)=(w(x),x+N),$ for each $x\in M.$ It is easy to verify that
$f$ is a monomorphism. It follows that $g(M)\leq g(N)+g(M/N).$ By induction on $\mu
(M),$ we complete the proof. \qed

\bigskip
If $R$ is a B\' ezout semi-CF-ring and $N$ an indecomposable cyclic $R$-module, then $J=\hbox{rad}(\hbox{ann}(N))=\{r\in R\mid\exists p\in \mathbb N,\hbox{ such that  }r^p\in \hbox{ann}(N)\}$ is the minimal prime ideal over
$\hbox{ann}(N).$ For every finitely generated module $M,$ from each
pure-composition series of $M,$ with indecomposable cyclic factors, we
can associate a sequence of prime ideals. Then we have the following
proposition:
\begin{proposition} 
\label{P:unilen}
Let $R$ be B\' ezout semi-CF-ring, and $M$
a finitely generated module. We consider two pure-composition
series of $M$ with indecomposable cyclic factors:

(s): $\{0\}=M_0\subset M_1\subset M_2\dots\dots\subset M_{n-1}\subset M_n = M.$

(s'): $\{0\}=M'_0\subset M'_1\subset M'_2\dots\dots\subset
M'_{n'-1}\subset M'_{n'} = M.$

We denote $A_i = \hbox{ann}(M_i/M_{i-1}),$ 
$A'_i = \hbox{ann}(M'_i/M'_{i-1}),$
$J_i=\hbox{rad}(A_i),$ $J'_i=\hbox{rad}(A'_i),$ $S=(J_i)_{1\leq i\leq n}$ and
$S'=(J'_i)_{1\leq i\leq n'}.$

Then $n=n',$ and there exists $\sigma $ in the symmetric group $S_n$
such that $J_{\sigma(i)}=J'_i$ for every $i, 1\leq i\leq n.$
\end{proposition}
\textbf{Proof.}  Let $C$ be a maximal totally ordered
subsequence of $S$ and $P$ a maximal ideal of $R$ containing every
term of $C$. For each $i$ such that $J_i\not\subseteq P$ we have
$(M_i)_P=(M_{i-1})_P,$ and for each $i$ such that $J_i\subseteq P$ we have
$(M_i)_P\supset(M_{i-1})_P.$ We deduce from (s) a pure-composition
series of $M_P,$ with cyclic factors, whose the length is equal to the
number of terms of $C$. By  \cite[Theorem 1.6, p. 177]{FuSa2}, the
pure-composition series of $M_P,$ deduced from (s) and (s') are
isomorphic, and this conclusion holds for each maximal ideal $P$
containing all terms of $C.$ Consequently there is a maximal totally ordered
subsequence $C'$ of $S',$ which has the same terms as $C.$
Conversely, the same conclusion holds for every maximal totally
ordered subsequence of $S'.$ 
Therefore $S$ and $S'$ have the same terms, eventually with
different indices. \qed

\bigskip
We will say that a sequence of ideals $(A_i)_{1\leq i\leq n}$ is
{\it almost totally ordered} if $A_i$ and $A_j$ are either comparable, or
comaximal, $\forall i,$ $\forall j,$ $1\leq i,j\leq n,$ and that this
sequence is {\it almost increasing} if either $A_i\subseteq A_j,$ or $A_i$
and $A_j$ are comaximal, $\forall i,$ $\forall j,$ $1\leq i<j\leq n,.$
Let us observe that every sequence of prime ideals of an arithmetic
ring is almost totally ordered.

 In the sequel, if $R$ is a B\' ezout semi-CF-ring and $M$ a finitely
 generated $R$-module, we denote $\ell (M)$ the length of every
 pure-composition series of $M,$ with indecomposable cyclic factors. 

\begin{proposition}
\label{P:semiso}
Let $R$ be a B\'ezout semi-CF-ring, and
$M$ a finitely generated $R$-module. Then the following assertions are true:
\begin{enumerate}
\item From every pure-composition series of $M,$ with indecomposable
  cyclic factors, we deduce a pure-composition series with an almost
  increasing annihilator sequence.
\item Any two pure-composition series of $M,$ with indecomposable
  cyclic factors and almost totally ordered annihilator sequences, are
  isomorphic.
\end{enumerate}
\end{proposition}
\textbf{Proof.}

(1) We consider the following pure-composition series of $M,$ with
  indecomposable cyclic factors,
\[\mathrm{(s):} 0=M_0\subset M_1\subset M_2\dots\dots\subset M_{n-1}\subset
M_n = M,\] where $M_i/M_{i-1}=R(x_i+M_{i-1}).$

We denote $A_i = \hbox{ann}(M_i/M_{i-1})$ and
$S$ the sequence of prime ideals associated to $M.$ We put
$S=(J_i)_{1\leq i\leq n}.$

We claim that, after a possible permutation of indices, we may assume
that $S$ is almost increasing. To prove this, we induct on $n.$ If
$n=1,$ it is obvious. Suppose $n>1.$ Let $C$ be a maximal totally
ordered subsequence of $S,$ and $S'$ the subsequence of all terms of
$S$ which are not terms of $C.$ If $C=S$ the claim is obvious. We may
assume that $C\not= S.$ Clearly, if $J_i$ is a term of $C$ and $J_k$ a
term of $S',$ then $J_i\subset J_k$ or $J_i$ and $J_k$ are
comaximal. Since the induction hypothesis can be applied to $S',$ then
it is possible to get an almost increasing sequence, if we begin to
index the terms of $C$ and we end with the terms of $S'.$

We will prove, by induction on $n=\ell (M),$ that we can get a
pure-composition series (t) of $M,$ with indecomposable cyclic factors:
$\{0\}\subset N_1\subset \dots N_{n-1}\subset M,$ with an almost
increasing annihilator sequence $(B_i)_{1\leq i\leq n}$ and such that
$\hbox{rad}(B_i)=J_i,$ for each $i,$ $1\leq i\leq n.$

If $n=1,$ this is obvious. We assume that $n>1.$ Let $k$ be the
smal\-lest index such that $\hbox{rad}(A_k)=J_n.$ If $k<n,$ one
of the five following cases holds:

(a): $\hbox{rad}(A_{k+1})$ and $J_n$ are comaximal. It follows that
$A_k$ and $A_{k+1}$ are also comaximal. Therefore, for every maximal
ideal $P,$ we have $(M_{k+1}/M_k)_P=\{0\}$ or $(M_k/M_{k-1})_P=\{0\}.$
Hence $M_{k+1}/M_{k-1}$ is cyclic and its annihilator is $A_k\cap
A_{k+1}.$ Then
$M_{k+1}/M_{k-1}\simeq R/(A_k\cap A_{k+1})\simeq R/A_k\oplus R/A_{k+1}.$
Hence $M_{k+1}/M_{k-1}=R(x'_k+M_{k-1})\oplus R(x'_{k+1}+M_{k-1}),$ where
$\hbox{ann}(x'_k+M_{k-1})=A_k$ and $\hbox{ann}(x'_{k+1}+M_{k-1})=A_{k+1}.$ Therefore $M_k$
can be replaced with $M'_k=M_{k-1}+Rx'_{k+1},$ where
$\hbox{ann}(M'_k/M_{k-1})=A_{k+1}$ and
$A_k=\hbox{ann}(M_{k+1}/M'_k).$ Then we replace $k$ with $(k+1)$ for the
next step.

(b): $\hbox{rad}(A_{k+1})\subset J_n.$ By applying \cite[Lemma 3.1]{ShWi} to $R/(A_k\cap A_{k+1}),$ we get that $A_{k+1}\subset A_k.$
With the same proof as in \cite[Lemma 1.3 p.172]{FuSa2}, we state that
$M_{k+1}/M_{k-1}=M_k/M_{k-1}\oplus R(x_{k+1}+M_{k-1}).$ Therefore $M_k$
can be replaced with $M'_k=M_{k-1}+Rx_{k+1},$ where
$\hbox{ann}(M'_k/M_{k-1})=A_{k+1}\subset
A_k=\hbox{ann}(M_{k+1}/M'_k).$ Then we replace $k$ with $(k+1)$ for the
next step.

(c): $\hbox{rad}(A_{k+1})=J_n$ and $A_k\subseteq A_{k+1}.$ Then we
replace $k$ with\ $(k+1)$ for the next step.

(d): $\hbox{rad}(A_{k+1})=J_n$ and $A_{k+1}\subset A_k.$ We do as in
the case (b).

(e): $\hbox{rad}(A_{k+1})=J_n$ and $A_{k+1}$ is not comparable with
$A_k.$ By applying Proposition~\ref{P:semlo} or~\ref{P:mipr} to
$M_{k+1}/M_{k-1},$ it follows that there exists a pure submodule
$M'_k$ of $M_{k+1},$ containing $M_{k-1},$ such that $M'_k/M_{k-1}$
and $M_{k+1}/M'_k$ are cyclic and $\hbox{ann}(M'_k/M_{k-1})=A_k\cap
A_{k+1}.$ For every maximal ideal $P$,
 $\{(A_k)_P,(A_{k+1})_P\}=\{(A_k\cap
A_{k+1})_P,(A_k+A_{k+1})_P\}.$ By using  \cite[Theorem 1.6,
p.177]{FuSa2}, it follows that the annihilator of $M_{k+1}/M'_k$ is 
$A_k+A_{k+1}.$ Therefore $M_k$ can be replaced with $M'_k.$ Then we
replace $k$ with $(k+1)$ for the next step.

After $(n-k)$ similar steps, we get from (s), a
pure-composition series (s') of $M,$ $\{0\}\subset M'_1\dots\subset
M'_{n-1}\subset M'_n = M,$ with the annihilator sequence
$(A'_i)_{1\leq i\leq n},$ and such that $\hbox{rad}(A'_n)=J_n.$ Let us
observe that, if $\hbox{rad}(A'_k)=J_n$, then $A'_k\subseteq A'_n,$
for each $k,$ $1\leq k<n.$ Hence, either $A'_i\subseteq A'_n,$ or
$A'_i$ and $A'_n$ are comaximal, for each $i,$ $1\leq i<n.$

Now $\ell (M'_{n-1})=n-1$ and $(J_i)_{1\leq i<n}$ is the sequence of
prime ideals associated to $M'_{n-1}.$ Moreover, (s') induces a
pure-composition series (s'') of $M'_{n-1},$ with the annihilator
sequence $(A'_i)_{1\leq i<n}.$ From the induction hypothesis, we get
from (s''), a pure-composition series (t') of $M'_{n-1},$ $\{0\}\subset
N_1\subset \dots N_{n-1}=M'_{n-1},$ with an almost increasing
annihilator sequence $(B_i)_{1\leq i<n}$ such that
$J_i=\hbox{rad}(B_i),$ for every $i,$ $1\leq i<n.$ We put $B_n=A'_n.$
From the above observation, it follows that $(B_i)_{1\leq i\leq n}$ is
an almost increasing sequence.

Let us observe that the case (e) is not possible if $(A_i)_{1\leq
  i\leq n}$ is almost totally ordered. In this case, (s) and (t) are
  isomorphic. 

(2) We consider two pure-composition series of $M,$ with indecomposable
cyclic factors and almost totally ordered annihilator sequences:

(s): $\{0\}=M_0\subset M_1\subset M_2\dots\dots\subset M_{n-1}\subset M_n = M.$

(s'):$\{0\}=M'_0\subset M'_1\subset M'_2\dots\dots\subset
M'_{n'-1}\subset M'_{n} =M.$ 

We denote $A_i = \hbox{ann}(M_i/M_{i-1})$ and
$A'_i = \hbox{ann}(M'_i/M'_{i-1}).$
From (1) we may assume that the annihilator sequences are almost
increasing and that $\hbox{rad}(A_i)=\hbox{rad}(A'_i)=J_i,$
$\forall i,$ $1\leq i\leq n.$  Let $P$ be a maximal ideal of $R.$ If
$J_i\subseteq P,$ then by \cite[Theorem 1.6, p.117]{FuSa2},
$(A_i)_P=(A'_i)_P.$ If $J_i\not\subseteq P,$ then $(A_i)_P=(A'_i)_P=R_P.$
Therefore $A_i=A'_i,$ for every $i,$ $1\leq i\leq n.$ \qed

\bigskip
Now, we can prove the implication (1) $\Rightarrow$ (2) of
Theorem~\ref{T:main}.
 
\textbf{Proof of Theorem~\ref{T:main}.}
(1) $\Rightarrow$ (2). By \cite[Corollary 3.9]{ShWi}, every
pure-composition series of a finitely generated $R$-module, with
indecomposable cyclic factors, has an almost totally ordered
annihilator sequence. Hence the result follows from
Proposition~\ref{P:pugo} and Proposition~\ref{P:semiso}. \qed

\bigskip
The following propositions are needed to prove the implication (2)
$\Rightarrow$ (1) of Theorem~\ref{T:main}. Now we suppose that $R$ is a
PCS-ring and if $M$ is a finitely generated module, we denote
$\ell(M)$ the length of every pure-composition series of $M$ with
indecomposable cyclic factors.

\begin{proposition} 
\label{P:gol}
Let $R$ be a PCS-ring. The following assertions
are true:
\begin{enumerate}
\item For each indecomposable cyclic module $N,\ g(N)=1$
\item For every finitely generated $R$-module $M,\ g(M)\leq\ell(M).$
\end{enumerate}
\end{proposition}
\textbf{Proof.}
\begin{enumerate}
\item Since $N\simeq R/\hbox{ann}(N),$ we can replace $R$ with
$R/\hbox{ann}(N)$ and assume that $R$ is an indecomposable module. Let
$I$ and $J$ be ideals of $R$ such that $I\cap J = \{0\}.$ Then
$R/I\oplus R/J$ is a finite direct sum of indecomposable cyclic
modules. If $I\not=\{0\}$ and $J\not=\{0\}$ then every annihilator of
each summand of $R/I\oplus R/J$ is a nonzero ideal. On the other hand we consider the following exact sequence~:

(i) $0\rightarrow R \overset{\sigma}{\rightarrow} R/I\oplus
R/J\overset{\varphi}{\rightarrow} R/(I+J) \rightarrow 0$ such that $\sigma(r) =
(r+I,r+J)$ and $\varphi(r+I,t+J) = (r-t)+(I+J).$ Let $P$ be a maximal ideal
of $R.$ Then we may assume that $I_P\subseteq J_P,$ and, if
$r(x+I_P,y+J_P) = (z+I_P,z+J_P)$, it is obvious that $r(x+I_P,x+J_P)
= (z+I_P,z+J_P).$ Hence $R_P$ is isomorphic to a pure $R_P$-submodule of
$(R/I\oplus R/J)_P.$ We deduce that (i) is a pure exact
sequence, and from a pure-composition series of $R/(I+J),$ we get a
pure-composition series of $R/I\oplus R/J$ such that the first annihilator
is $\{0\}.$ Consequently, if $I\not=\{0\}$ and $J\not=\{0\},$ $R/I\oplus
R/J$ has two pure-composition series which are not isomorphic. Since $R$ is
a PCS-ring, we deduce that $R$ is a uniform module and that $g(R) = 1.$

\item By induction on $n=\ell(M).$ If $n=1$ the result follows from
  (1). Now suppose $n>1.$ Let $\{0\}\subset M_1\subset\dots\subset
M_{n-1}\subset M_n = M$ be a pure-composition series of $M,$ with
indecomposable cyclic factors. As in the proof of Proposition~\ref{P:pugo}, we
get that $g(M)\leq g(M_{n-1})+1,$ because $g(M/M_{n-1})=1.$ Since $\ell(M_{n-1})=n-1,$ it follows
from the induction hypothesis that $g(M)\leq n.$ \qed
\end{enumerate}
\begin{proposition} 
\label{P:pri}
Let $R$ be a PCS-ring. Then $R$ has only finitely
many minimal prime ideals.
\end{proposition}
\textbf{Proof.}
We may assume that $R$ is a reduced ring. If $S$ denotes the total ring of
quotients of $R,$ then by \cite[Proposition 2, p. 106]{Lam}, $S$ is a
Von Neumann regular ring. By Proposition~\ref{P:gol}, $g(R)$ is
finite. We deduce from \cite[Proposition 2, p. 103]{Lam} that $S$ is
semi-local. Consequently $S$ is a semi-simple ring, i.e. $S =
\prod_{i=1}^n K_i,$ where $K_i$ is a field for every $i,\
1\leq i\leq n.$ If $u : R\rightarrow S$ is the natural map, and,
$\forall i,\ 1\leq i\leq n,$ $\ p_i : S\rightarrow K_i$ the canonical
epimorphism, we denote $P_i = \ker(p_i\circ u).$ Then $P_i$ is a prime
ideal and $\bigcap_{i=1}^n P_i = \{0\}.$ We deduce that
$\{P_i\mid 1\leq i\leq n\}$ is the set of minimal prime ideals of $R.$ \qed
\begin{proposition} 
\label{P:puni}
Let $R$ be a PCS-ring with a unique minimal prime
ideal $Q.$ Then $Q$ is a uniserial module.
\end{proposition}
\textbf{Proof.}
If there exist $a$ and $b$ in $Q$ such that $a\notin Rb$ and $b\notin Ra,$
then $R/(Ra\cap Rb)$ is indecomposable with $g(R/(Ra\cap Rb))\geq 2.$ From
Proposition~\ref{P:gol}, we get a contradiction. \qed

\begin{proposition} 
\label{P:onema}
Let $R$ be a PCS-ring. Then every nonminimal
prime ideal is contained in only one maximal ideal.
\end{proposition}
\textbf{Proof.} We do a similar proof as in \cite[Lemma 10]{WiWi}. By
Proposition~\ref{P:pri}, $R$ has only finitely many minimal prime ideals. Since
$R$ is arithmetic, these minimal prime ideals are pairwise comaximal
and consequently $R$ is a finite product of indecomposable PCS-rings
with a unique minimal prime. Hence we may assume that $R$ has only one
minimal prime ideal $Q.$ Let $P$ be a prime ideal such that $Q\subset
P$ and $a\in P\setminus Q.$ Then, from \cite[Lemma 3.1]{ShWi}, it
follows that $Q\subset Ra.$
We may also assume that $P$ is minimal over $Ra.$ Suppose that $P\subseteq
M_1\cap M_2$ where $M_1$ and $M_2$ are two distinct ma\-ximal ideals of $R.$
Let $x_1$ and $x_2$ in $R,\ x_i\in M_j\setminus M_i,\ i\not= j,$ such that
$x_1+x_2 = 1.$ If $S=R\setminus(M_1\cup M_2),$ we denote $R' = S^{-1}R.$
Then $ax_i\notin aM_iR',\ \forall i,\ i=1$ or 2. Else, $ax_i =
\displaystyle{\frac{am_i}{s}}$ where $m_i\in M_i$ and $s\in S.$ Hence there
exists $t\in S$ such that $a(tsx_i-tm_i) = 0.$ It follows that
$(tsx_i-tm_i)\in Q,$ and since $tm_i\in M_i,$ we get that $tsx_i\in M_i.$
But $t,s$ and $x_i$ are not in $M_i,$ hence we get a contradiction. Let $I
= aM_1\cap aM_2.$ Then $\displaystyle{\frac{ax_i}{1}}\notin R'I.$ Now we
consider $R_1 = R/I.$ Then $R_1$ is also a finite product of indecomposable
rings with a unique minimal prime ideal. Hence we may assume that $P_1 =
P/I$ is the only minimal prime ideal of $R_1.$ Then, if $S_1 =
R_1\setminus(M_1/I\cup M_2/I)$ we have $S_1^{-1}R_1\simeq R'/R'I.$
Consequently $ax_i+I \not= I,\ \forall i,\ i\in\{1,2\}.$ Since $R$ is a
B\'ezout ring and $x_1+x_2 = 1,\ Rax_1\cap Rax_2 =  Rax_1x_2.$ We deduce
that $R_1(ax_1+I)\cap R_1(ax_2+I) =\{0\}.$ This contradicts that $P_1$ is
uniserial by Proposition~\ref{P:puni}. \qed

\bigskip
Now, we complete the proof of Theorem~\ref{T:main}.

\textbf{Proof of Theorem~\ref{T:main}.}
(2) $\Rightarrow$ (1). By Proposition~\ref{P:Kap}, $R$ is a B\'ezout
ring. We deduce from Propositions~\ref{P:pri}, \ref{P:puni}
and~\ref{P:onema} and from \cite[Theorem 3.12]{ShWi} that $R$ is also a
CF-ring. \qed

\bigskip
It is easy to get a semilocal arithmetic ring $R$ which is not a
PCS-ring.

\begin{example} 
\label{E:npcs1}
\textnormal{Let $D$ a semilocal and h-local B\'ezout domain, $F$
its classical field of quotients, $P$ and $Q$ two distinct maximal ideals of
$D$,
 $E =(F/PD_P)\oplus( F/QD_Q)$ and $R =
\displaystyle\Bigl\{{\binom{d\ \ e}{0\ \ d}}\mid d\in D, e\in
E\Bigr\}$ the trivial extension of $D$ by $E.$ Then $R$ is a semilocal
arithmetic ring with a unique minimal prime ideal which is not uniserial.
By Proposition~\ref{P:puni}, $R$ is not a PCS-ring.}
\end{example}
\begin{example} 
\label{E:npcs2}
\textnormal{Let $R$ be the Kaplansky domain of \cite[example 1]{Hen}
. We may assume that $R$ is semilocal and nonlocal. Then
its Jacobson radical is a nonzero prime ideal contained in every
maximal ideal. Since $R$ is not h-local, $R$ is not a PCS-ring.}
\end{example}
The following questions can be proposed~:

Let $R$ be a ring such that every finitely generated $R$-module has a 
pure-composition series with indecomposable cyclics factors.
\begin{enumerate}
\item Is $g(M)$ finite, for every finitely generated $R$-module $M?$
\item Is $Minspec(R)$ a finite set ? $Minspec(R)$
denotes the set of minimal prime ideals of $R.$
\item For every finitely generated $R$-module $M,$ are the lengths of
  any two pure-composition series of $M,$ with indecomposable cyclic
  factors, equal?
\end{enumerate}

Let us observe that, if the answer to the question (1) is positive,
then $Minspec(R)$ is finite (see the proof of Proposition~\ref{P:pri}), and if
$Minspec(R/I)$ is finite, for every proper ideal $I,$ then the
question (3) also has a positive answer (see the proof of
Proposition~\ref{P:unilen}.)

When $Minspec(R)$ is compact in its Zariski topology, we can answer
to the question (2).

\begin{proposition} 
\label{P:arifin}
Let $R$ be an arithmetic ring. We suppose that
$Minspec(R)$ is compact and every finitely generated module
has a pure-composition series with indecomposable cyclic factors.

Then $Minspec(R)$ is a finite set.
\end{proposition}
\textbf{Proof.}
Since $R$ has a pure-composition series with indecomposable cyclic factors, $R$ is a finite direct product of indecomposable Kaplansky rings. We
may assume that $R$ is semi-prime and indecomposable. Since $R_P$ is a
valuation domain for every maximal ideal $P,$ by \cite[Proposition 10]{Que},
 $R$ is semi-hereditary. Let $r\in R,\ r\not= 0.$
 Since $\hbox{ann}(r)$ is generated by an idempotent,
 $\hbox{ann}(r)=\{0\}.$ Hence $R$ is a domain. \qed

\bigskip
We can answer to these questions when $R$ is a Von Neumann regular ring.
Then every exact short sequence of $R$-modules is pure. It follows that
every finitely generated $R$-module has a pure-composition series with
cyclic factors. However, we have the following result.
\begin{proposition} 
\label{P:VN}
Let $R$ be a Von Neumann regular ring. Then the
following statements are true.
\begin{enumerate}
\item A finitely generated $R$-module has a pure-composition series
with indecomposable cyclic factors if and only if it is a semi-simple module.
\item Every finitely generated $R$-module has a pure-composition
series with indecomposable cyclic factors if and only if $R$ is a
semi-simple ring. 
\end{enumerate}
\end{proposition}
\textbf{Proof.}
Clearly (1) $\Rightarrow$ (2). We can also deduce (2) from
Proposition~\ref{P:arifin}. (1) is an immediate consequence of the
following Theorem. \qed 
\begin{theorem} 
\label{T:VN}
Let $R$ be a ring. Then the following assertions are
equivalent~:
\begin{enumerate}
\item $R$ is a Von Neumann regular ring.
\item Every indecomposable $R$-module is simple.
\end{enumerate}
\end{theorem}
\textbf{Proof.}
(2) $\Rightarrow$ (1).  For every simple $R$-module $S,\ E_R(S)\simeq
S.$ Since every simple $R$-module is injective, we deduce that $R$ is Von Neumann regular.

(1) $\Rightarrow$ (2). Let $M$ be an indecomposable $R$-module, $M\not=\{0\}.$
Then there exists a maximal ideal $P$ such that $M_P\not=\{0\}.$ Let $a\in
P.$ Then there exists an idempotent $e$ such that $Ra=Re.$ Since
$M_P\not=\{0\}$ and $(1-e)\notin P,$ we have $(1-e) M\not=\{0\}.$ We deduce
that $eM = \{0\}$ since $M$ is indecomposable. Consequently $M$ is an
$R/P$-module. Since $R/P$ is a field, $M$ is a simple module. \qed

\bigskip
By using this theorem, we can also answer to the following question proposed
by R. Wiegand in \cite{Wie}~: are the following assertions equivalent for a
commutative ring $R$ and an integer $n\geq 1$ ?
\begin{enumerate}
\item Every finitely generated $R$-module is a finite direct sum of submodules
generated by at most $n$ elements.
\item $\mu(M)\leq n,$ for every indecomposable finitely generated $R$-module
$M.$
\end{enumerate}

If $R$ is a Von Neumann regular ring, then $R$ satisfies the second
assertion by Theorem~\ref{T:VN}, with $n=1.$ But if $R$ satisfies the first
assertion, then $R$ is semi-simple by \cite[Corollary 21.7]{Pie} . Consequently
the answer to Wiegand's question is negative.

\section{The Goldie dimension} 
\label{S:gol}
In this section we will prove the following theorem.

\begin{theorem} 
\label{T:cyc}
Let $R$ be a PCS-ring. Then a finitely generated
$R$-module $M$ is a direct sum of cyclic submodules if and only if $g(M) =
\ell(M).$
\end{theorem}

When $R$ is a valuation ring this theorem was proved by L. Salce and P.
Zanardo \cite[Corollary 3.5]{SaZa}, and when $R$ is an h-local B\'ezout domain,
it was proved by C. Naud\'e \cite[Theorem 2.2]{Nau}. Some preliminary
results are needed to prove this theorem.

\begin{lemma} 
\label{L:ann} 
Let $R$ be an arithmetic ring, $\{0\}\subset
M_1\subset\dots\subset M_{n-1}\subset M_n = M$ a pure-composition series of
a finitely generated $R$-module $M,$ with $M_i/M_{i-1}
= R(x_i+M_{i-1})$ for all $i,\ 1\leq i\leq n,$ and with an increasing
annihilator sequence $(A_i)_{1\leq i\leq n}.$ Let $c_1,\dots,c_n$ in
$R$ such that $\sum_{i=1}^n c_ix_i = 0.$

 Then $c_i\in A_n,$ for every $i,\ 1\leq i\leq n.$
\end{lemma}
\textbf{Proof.}
We prove this lemma by induction on $n.$ If $n=1$ it is obvious. If $n>1,$
then $c_n\in A_n.$ Moreover, since $M_{n-1}$ is a pure submodule of  $M,$
there exist $d_1,\dots,d_{n-1}$ in $R$ such that $c_nx_n =
c_n(\sum_{i=1}^{n-1}d_ix_i) = -\sum_{i=1}^{n-1}
c_ix_i.$ We deduce the following equa\-lity
$\sum_{i=1}^{n-1}(c_nd_i+c_i)x_i =~0.$ From the induction
hypothesis it follows that $(c_nd_i+c_i)\in A_{n-1}$ for every $i,\ 1\leq
i\leq n-1.$ Since $A_{n-1}\subseteq A_n$ and  $c_n\in A_n,$ we deduce that
$c_i\in A_n,\ \forall i,\ 1\leq i\leq n.$ \qed

\begin{lemma} 
\label{L:gol1}
Let $R$ be a B\'ezout and a torch ring, $Q$ its minimal
prime ideal and $P$ the only maximal ideal such that $Q_P\not=\{0\}.$ Let
$M$ be a finitely generated $R$-module, $\{0\}\subset M_1\subset
\dots\subset M_n = M$ a pure-composition series, with 
$M_i/M_{i-1}= R(x_i+M_{i-1})$ for all $i,$ $ 1\leq i\leq
n,$ and with an increasing annihilator sequence $(A_i)_{1\leq i\leq n}.$
We assume that $A_n\subseteq Q.$ Let $M' = M_{n-1}.$

If $M'$ is not essential in $M,$ then $M$ has a non-zero cyclic summand.
\end{lemma}
\textbf{Proof.}
Let $y\in M\setminus M'$ such that $Ry\cap M' = \{0\}.$ Then $y =
\sum_{i=1}^n a_ix_i,$ where $a_i\in R$ for all $i,\ 1\leq
i\leq n.$ By Propositions~\ref{P:mipr} and~\ref{P:Kap} 
$R$ is a Kaplansky ring, hence there exist $a, b_1,\dots,b_n$ in
$R$ such that $a_i=ab_i$ for every $i,\ 1\leq i\leq n,$ and
$\sum_{i=1}^n Rb_i = R.$ We put $z = \sum_{i=1}^n
b_ix_i.$

We claim that $M/QM$ is a free $R/Q$-module with basis
$\{x_i+QM\mid 1\leq i\leq n\}.$ Let $c_1,\dots,c_n\in R$ such that
$\sum_{i=1}^n c_ix_i\in QM.$ There exist $q_1,\dots,q_n$ in
$Q$ such that $\sum_{i=1}^n (c_i-q_i)x_i = 0.$ From
Lemma~\ref{L:ann}, since $A_n\subseteq Q,$ it follows that $c_i\in Q,$
 $\forall i,1\leq i\leq n.$  

Now we claim that $R_Py\not=\{0\}.$ Else there exists $s\notin P$ such that
$sy=0.$ From Lemma~\ref{L:ann}, it follows that $a_i\in Q,\ \forall i, \ 1\leq
i\leq n,$ whence $y\in QM.$ If $P'$ is a maximal ideal, $P'\not= P,$ then
$M_{P'} \simeq (M/QM)_{P'}$ and consequently
$R_{P'}y = \{0\}.$ Since $Ry\not=\{0\},$ we get a contradiction.

Since $R_Py\not=\{0\},\ R_Py\cap M'_P = \{0\}$ and $R_P$ is a valuation
ring, it follows that $R_Pz\cap M'_P = \{0\}.$

Since $R$ is a Kaplansky ring, there exists an invertible $n\times n$ matrix
$\Lambda$ such that $[b_1,\dots,b_n]\Lambda=[1,0,\dots,0].$

Then $z=
[b_1,\dots,b_n][x_1,\dots,x_n]^t=
[b_1,\dots,b_n]\Lambda\Lambda^{-1}
[x_1,\dots,x_n]^t.$

 Hence $z=[1,0,\dots,0]\Lambda^{-1} [x_1, \dots,x_n]^t.$ 
We put $[z_1,\dots,z_n]^t=\Lambda^{-1} [x_1,\dots,x_n]^t.$ Thus $z=z_1$ and
$\{z_1,\dots,z_n\}$ generates $M.$ Let $N$ be the submodule of $M$
generated by  $\{z_2,\dots,z_n\}.$ Then $\{z_i+QM\mid 1\leq i\leq n\}$ is a
basis of $M/QM,$ and for every maximal ideal $P', P\not= P',$ we have
$M_{P'} = (Rz_1)_{P'}\oplus N_{P'}.$

We put $\Lambda = (\lambda_{i,k})_{1\leq i,k\leq n}.$ Then,
since $R_P$ is a valuation ring, there exists $j,$ $1\leq j\leq n,$ such
that $\lambda_{j,1}$ and $b_j$ are units of $R_P$ and such that $\Lambda_{j,1}$ is an
invertible matrix, where $\Lambda_{j,1}$ is the $(n-1)\times(n-1)$ matrix
obtained from $\Lambda,$ by deleting its first column and its $j$-th
row.

Moreover $\Lambda_{j,1} [z_2,\dots,z_n]^t=$
\[[x_1-\lambda_{1,1}z_1,\dots,x_{j-1}-\lambda_{j-1,1}z_1,x_{j+1}-\lambda_{j+1,1}z_1,\dots,x_n-\lambda_{n,1}z_1]^t\]
 and $\{x_i-\lambda_{i,1}z_1\mid 1\leq i\leq n,\ i\not= j\}$ also generates $N_P.$ Since $b_j$ is a unit of $R_P,$
as in the proof of  of \cite[Lemma 2.1, p. 179]{FuSa2}, we get that 
$R_Pz_1\cap\Bigl(\displaystyle\sum_{\scriptstyle
    i=1\atop\scriptstyle i\not=j}^n R_Px_i\Bigr) = \{0\}.$ 
Let $x\in (Rz_1)_P\cap N_P.$ Then there exist $c$
and $c_i$ in
$R_P,$ for every $i,\ 1\leq i\leq n,\ i\not= j,$ such that
$x=cz_1=\displaystyle\sum_{\scriptstyle i=1\scriptstyle \atop\scriptstyle
i\not= j}^n c_i(x_i-\lambda_{i,1}z_1).$ From above, we deduce that
$\Bigl(c+\displaystyle\sum_{\scriptstyle i=1\atop\scriptstyle i\not= j}^n
c_i\lambda_{i,1}\Bigr)z_1 = \sum_{{\scriptstyle i=1\atop\scriptstyle i\not=
j}}^n c_ix_i = 0.$ From Lemma~\ref{L:ann}, it follows that $c_i\in(A_n)_P$ for
every $i,\ 1\leq i\leq n,\ i\not=j.$ On the other hand, if we set $d =
c+\displaystyle\sum_{\scriptstyle i=1\atop\scriptstyle i\not= j}^n
c_i\lambda_{i,1},$ then $dz_1=0$ implies that $db_i\in(A_n)_P,\ \forall i,\
1\leq i\leq n.$ Since $b_j$ is a unit of $R_P,$ it follows that
$d\in(A_n)_P,$ whence $c\in (A_n)_P.$ Consequently $cz_1\in M'_P.$ Since
$M'_P\cap R_Pz_1 = \{0\},$ we also get that  $R_Pz_1\oplus N_P = M_P.$ \qed

\begin{lemma} 
\label{L:gol2}
Let $R$ be a B\'ezout and a torch ring and $Q$ its minimal
prime ideal. Let $M$ be a finitely generated module and  $\{0\}\subset
M_1\subset M_2\subset\dots\subset
M_n = M$ a pure-composition series with  $M_i/M_{i-1}=
R(x_i+M_{i-1})$ for all $i,\ 1\leq i\leq n,$ and with an increasing
annihilator sequence $(A_i)_{1\leq i\leq n}.$ We assume that $Q\subset
A_n.$ If  $M/M_{n-1}= R(x'_n+M_{n-1})\oplus
R(x"_n+M_{n-1}),$ where $R(x'_n+M_{n-1})$ is indecomposable, we set $M' =
M_{n-1}+Rx"_n.$

If $M'$ is not essential in $M,$ then $M$ has a nonzero cyclic summand.
\end{lemma}
\textbf{Proof.}
There exists $y\in M\setminus M'$ such that $M'\cap Ry = \{0\}.$ Let $A' =
\hbox{ann}(M/M').$ Then $A_n\subseteq A'$ and there exists only one
maximal ideal $P$ such that $A'\subseteq P.$ Therefore, if $P'$ is a
maximal ideal, $P'\not= P,$ we have $(M/M')_{P'} = \{0\}.$ It follows
that there exists $s\notin P'$ such that $sx'_n\in M'.$ We also get 
that $sy\in M'$ and since $M'\cap Ry = \{0\},$ we have $sy=0.$
Consequently $R_{P'}y = \{0\}$ and $M'_{P'} = M_{P'}$ for every maximal ideal $P'\not= P.$

There exist $a_1,\dots,a_n$ in $R$ such that $y = \sum_{i=1}^n
a_ix_i.$ As in the previous lemma, there exist $a, b_1,\dots,b_n$ such that
$a_i = ab_i$ for every $i,\ 1\leq i\leq n,$ and such that
$\sum_{i=1}^n Rb_i = R.$ Then $a\notin A',$ else $y\in M'.$
 Since $a\notin Q$ and $R/Q$ is an $h$-local B\'ezout
domain, it follows that $R/Ra = \bigoplus_{j=1}^m R/Rc_j,$
where $R/Rc_j$ is indecomposable for every $j, \ 1\leq j\leq m.$ We
denote $P_j$ the only maximal ideal of $R$ such that $Rc_j\subseteq P_j,$
for each $j,\ 1\leq j\leq m.$

If $P\not= P_j$ for each $j,\ 1\leq j\leq m,$ then $a$ is a unit of $R_P.$
In this case we set $z=y$ and $d_i=a_i,$ for every $i,\ 1\leq i\leq n.$

If there exists $j,\ 1\leq j\leq m,$ such that $P=P_j,$ then we put
$c=c_j.$ There exists $d\in R$ such that $a=cd,$ and $d$ is a unit of
$R_P.$ We set $d_i=b_id,$ for every $i,\ 1\leq i\leq n,$ and we put $z =
\sum_{i=1}^n d_ix_i.$ For every maximal ideal $P'\not= P,\ c$
is a unit of $R_{P'},$ whence $R_{P'}z = R_{P'}y = \{0\}.$ Since $M'_P\cap
R_Py = \{0\}$ and $R_P$ is a valuation ring, it follows that
$R_Pz\cap M'_P = \{0\}.$

In the two cases, there exists $j,\ 1\leq j\leq n,$ such that $d_j$ is a
unit of $R_P.$ Let $\ell$ be the greatest index such that $A_\ell\subseteq Q.$
We claim that $\ell<j.$ Else there exist $t\in PR_P,$ and $d'_{\ell+1},\dots,d'_n\in R_P$
such that $d_i = td'_i$ for every $i,\ \ell<i\leq n.$ Let $s\in
R_PA'\setminus R_PQ.$ Then $sz=0.$ We get the following equality~:
$st(\sum_{i=\ell+1}^n d'_ix_i) =
-s(\sum_{i=1}^\ell d_ix_i).$ Since $(M_\ell)_P$ is a
pure-submodule of $M_P,$ there exist $d'_1,\dots,d'_\ell$ in $R_P$ such
that $\sum_{i=1}^\ell s(td'_i+d_i)x_i = 0.$ From
Lemma~\ref{L:ann} it follows that $s(td'_i+d_i)\in R_PA_\ell\subseteq
R_PQ,$ for every $i,\ 1\leq i\leq\ell.$ Since $s\notin R_PQ$ and $t\in
PR_P,$ we deduce that $d_i\in PR_P$ for every $i,\ 1\leq i\leq\ell.$
Hence we get a contradiction.

Since $Q\subset A_j\subseteq A_n,$ then  $M_j/M_{j-1}=
R(x'_j+M_{j-1})\oplus R(x"_j+M_{j-1}),$ such that 
$(M_j/M_{j-1})_P =
R_P(x'_j+(M_{j-1})_P),$ $R_P(x"_j+(M_{j-1})_P) = \{0\}$ and  $R_{P'}(x'_j+(M_{j-1})_{P'}) =
\{0\}$ for every  maximal ideal $P'\not= P.$

Let $N = \displaystyle\sum_{\scriptstyle i=1\atop \scriptstyle i\not=j}^n
Rx_i+Rx"_j.$ Then $N_P = \displaystyle\sum_{\scriptstyle i=1\atop
\scriptstyle i\not= j}^n R_P x_i$ and $N_{P'} = M_{P'}$ for every maximal
ideal $P'\not= P.$ As in the proof of \cite[Lemma 2.1, p. 179]{FuSa2}, we state
that $N_P\oplus R_Pz = M_P.$ Consequently we get $M = N\oplus Rz.$ \qed

\begin{proposition} 
\label{P:cyc}
Let $R$ be a B\'ezout and a torch ring. Then  every
finitely $R$-module $M$ contains a pure and essential submodule $B$ which
is a direct sum of $g(M)$ cyclic modules.
\end{proposition}
\textbf{Proof.}
We induct on $m = \ell(M).$ The case $m=1$ is obvious. Let $m>1$ and
$\{0\}\subset M_1\dots\subset M_n = M$ be a pure-composition series with cyclic
factors and an increasing annihilator sequence $(A_i)_{1\leq i\leq n}.$
Let $Q$ be the minimal prime of $R.$ If $A_n\subseteq Q,$ then $n=m$ and we
set $M' = M_{n-1}.$ If $Q\subset A_n,$ then  $M/M_{n-1}
= R(x'_n+M_{n-1})\oplus R(x"_n+M_{n-1})$ where $R(x'_n+M_{n-1})$  is
indecomposable ; in this case we set $M' = M_{n-1}+Rx"_n.$ In the two
cases, $\ell(M') = m-1.$~

If $M'$ is essential, $M'$ has a pure and essential submodule $B'$ which is
a direct sum of $g(M')$ non-zero indecomposable cyclic  modules ; in this
case, we are done by setting $B = B'.$

If $M'$ is not essential, then by Lemma~\ref{L:gol1} or Lemma~\ref{L:gol2},  $M = N\oplus
Rz,$ for some $z\in M,\ z\not=0,$ where $Rz$ is indecomposable and $N$ a
submodule of $M.$ Then $\ell(N) = m-1.$ Thus $N$ has a pure and essential
submodule $B"$ which is a direct sum of $g(N)$ non-zero indecomposable
cyclic submodules. We put $B = B"\oplus Rz$ to conclude the proof. \qed

\bigskip
Now we can prove the Theorem~\ref{T:cyc}.

\textbf{Proof of Theorem~\ref{T:cyc}.}
By Proposition~\ref{P:gol}, the Goldie dimension of every cyclic indecomposable
module is one. Consequently if $M$ is a finite direct sum of cyclic
modules, we have $\ell(M) = g(M).$

Conversely let $M$ be a finitely generated module such that $\ell(M) =
g(M).$ By Proposition~\ref{P:pri}, $R = \prod_{j=1}^m R_j,$ where
$R_j$ is an indecomposable PCS-ring for every $j,\ 1\leq j\leq m.$ We
deduce that  $M\simeq \prod_{j=1}^m M_j,$ where $M_j
= R_j\otimes_R M.$ Then $\ell(M) = \sum_{j=1}^m\ell(M_j)$  and
$g(M) = \sum_{j=1}^m g(M_j).$ From Proposition~\ref{P:gol}, we deduce
that $\ell(M_j) = g(M_j)$ for every $j,\ 1\leq j\leq m.$ Consequently we
may assume that $R$ is indecomposable. When $R$ is a valuation ring, the
result was proved by L. Salce and P. Zanardo \cite[Corollary
3.5]{SaZa}, (or \cite[Theorem 2.4, p. 180]{FuSa2}), and when $R$ is an $h$-local B\'ezout domain, it was
proved by C. Naud\'e \cite[Theorem 2.2.]{Nau}. Hence we may assume
that $R$ is a B\'ezout and torch ring. We prove the result by using
Proposition~\ref{P:cyc} and Proposition~\ref{P:gol}, with the same proof as in 
\cite[Theorem 2.4, p. 180]{FuSa2}. \qed

\bigskip
Let $R$ be a B\' ezout semi-CF-ring, $M$ a finitely generated $R$-module and 
$(s):$ 
$\{0\}=M_0\subset M_1\subset \dots \subset M_{n-1}\subset M_n=M,$ a
pure-composition series of $M,$ with indecomposable cyclic factors. We
put $g((s))=\sum_{i=1}^n g(M_i/M_{i-1}).$ As in the proof
of Proposition~\ref{P:gol}, we state that $g(M)\leq g((s)).$

However, if $(s')$ is another pure-composition series of $M,$ with
indecomposable cyclic factors, we have not necessarely
$g((s))=g((s')).$ For instance:
\begin{example} 
\textnormal{Let $R$ be the ring defined in
    example~\ref{E:npcs1}. We put $E_1=F/PD_P,$ $E_2=F/QD_Q,$
    $A_i=\displaystyle\Bigl\{\binom{0\ \ e}{0\ \ 0}\mid e\in
    E_i\Bigr\},$ where $i=1$ or $2,$ and $J=A_1+A_2.$ Then $A_1\cap
    A_2=\{0\}$ and $J$ is the minimal prime ideal of $R.$ Let
    $M=R/A_1\oplus R/A_2$ and $(s)$ be the pure-composition series
    with the annihilator sequence $(A_1,A_2).$ Then  $g(M)=g((s))=2.$
    But $M$ has a pure-composition series $(s')$ whose factors are $R$
    and $R/J.$ We have $g((s'))=3.$ Let us observe that $M$ is not
    isomorphic to $R\oplus R/J.$}
\end{example}

If $M$ is a finitely generated module over a B\' ezout semi-CF-ring $R,$ we
denote $\mathcal S (M)$ the set of all pure-composition series of $M,$
with indecomposable cyclic factors, and we put $h(M)=inf\{g((s))\mid
(s)\in \mathcal S (M)\}.$ Then we have the following proposition:
\begin{proposition} Let $R$ be a B\' ezout semi-CF-ring and $M$ a 
finitely generated $R$-module. Then the following assertions are true:
\begin{enumerate}
\item $\mu (M)\leq\ell (M)\leq h(M).$
\item $g(M)\leq h(M).$
\item If $M$ is a direct sum of cyclic modules, then $g(M)=h(M).$
\end{enumerate}
\end{proposition}

We don't know if the converse of the third assertion holds.

\section{Pure-injectivity and RD-injectivity}
\label{S:inj}

First, for every integer $n\geq 2,$ we give an example of an 
artinian module $M$ over a noetherian domain $R,$ which has an 
RD-composition series of length $n,$ with uniserial factors. 

\begin{example} 
\label{E:art}
\textnormal{Let $K$ be an algebraically closed field of
characteristic 0, $K[X,Y]$ the polynomial ring in two variables $X$ and
$Y,$ and $f(X,Y) = Y^n-X^n(1+X),$ where $n\in\mathbb N,\ n\geq 2.$ By
considering that $f(X,Y)$ is a polynomial in one variable $Y$ with
coefficients in $K[X],$ it follows from Eisenstein's criterion that $f(X,Y)$ is
irreducible. Then $R = \displaystyle{\frac{K[X,Y]}{f(X,Y)K[X,Y]}}$ is a
domain. Let $x$ and $y$ be the  images of $X$ and $Y$ in $R$ by the natural
map and $P$ the maximal ideal of $R$ generated by $\{x,y\}.$ If $\widehat
R$ is the completion of $R$ in its $P$-adic  topology, then $\widehat
R\simeq\displaystyle{\frac{K[[X,Y]]}{f(X,Y)K[[X,Y]]}}.$ Since $K$ has
$n$ distinct $n$-th roots of unity, by applying Hensel's Lemma to
$K[[X]],$ we deduce that there exist $u_1(X),\dots,u_n(X)$ in $K[[X]]$
such that $f(X,Y) =\prod_{i=1}^n(Y-Xu_i(X)).$ Let $E = E_R(R/P).$ Then by
\cite{Matl}, $\hbox{End}_RE = \widehat R,\ E$ is also an injective
$\widehat R$-module and every $R$-submodule of $E$ is also an
$\widehat R$-submodule. Moreover, there is a bijection between the set
of ideals of $\widehat R$ and the set of submodules of $E.$
For every $k,$ $1\leq k\leq n,\ Q_k = (y-xu_k(x))\widehat R$ is a
minimal prime ideal of $\widehat R,$ and $Q_k\cap R = \{0\}.$}

\textnormal{For every $k,\ 1\leq k\leq n,$ let $F_k = \{e\in E\mid
Q_k\subseteq\hbox{ann}_{\widehat R}(e)\}.$ Then $\widehat
R/Q_k \simeq K[[Z]],$ and consequently it is a discrete valuation
domain. We deduce that $F_k$ is a uniserial $R$-module. Since $F_k$ is an
injective $\widehat R/Q_k$-module, it follows that
$F_k$ is a divisible $R$-module for every $k,\ 1\leq k\leq n.$ We set $E_k
= \sum_{i=1}^k F_i$ for every $k,\ 1\leq k\leq n.$ Then
$\hbox{ann}_{\widehat R}(E_k) = \bigcap_{i=1}^k Q_i$ and
consequently
$E_n = E.$ Moreover,  $E_k/E_{k-1}=
(E_{k-1}+F_k)/E_{k-1}
\simeq F_k/(F_k\cap E_{k-1}).$ We deduce that 
$E_k/E_{k-1}$ is a uniserial $R$-module, and
since $R$ is a domain it follows that
$E_k/E_{k-1}$ is also divisible. Consequently
$E_k$ is divisible over $R,$ $\forall k,$ $1\leq k\leq n,$ and the
following chain  $\{0\}\subset E_1\subset\dots\subset E_{n-1}\subset E_n=E,$ is
an RD-composition series of $E,$ of length $n,$ with uniserial
factors.}
\end{example}

\bigskip
Now, we give some examples of pure-injective modules that fail to be
RD-injective, over noetherian domains.

\begin{example} 
\textnormal{Suppose that $R,\ F_k$ and $E_k,\ 1\leq k\leq n,$
are defined as in the example~\ref{E:art}. Then, since $F_k$ and $E_k$ are
artinian $R$-modules, they are also pure-injective by
\cite[Proposition 9]{War1}. Since $E$ is indecomposable, $F_k$ and
$E_k$ $(k\leq n-1)$ are not RD-injective.}
\end{example}

\bigskip
Recall that an $R$-module $M$ is {\it linearly compact} (in its discrete
topology) if for every family of cosets $\{x_i+M_i\mid i\in I\},$ with the
finite intersection property, has a non-void intersection. 

We can also find a noetherian and linearly compact module $M$ which is not
RD-injective over a noetherian domain $R.$

\begin{example} \textnormal{Let $R$ be a complete local regular
    noetherian domain $R$ with Krull dimension $n\geq 2$ (for example,  $R = K[[X_1,\dots,X_n]],$
where $K$ is a field). Then there exists an exact sequence:
\[0 \longrightarrow F_n\overset{u_n}{\longrightarrow} F_{n-1}\overset
{u_{n-1}}{\longrightarrow} F_{n-2}\dots \dots F_1 \overset
{u_1}{\longrightarrow} F_0,\] 
where $F_i$ is a free module of finite rank, for every $i,\ 0\leq i\leq n,$
and such that Im $u_{n-1}$ is not a projective module. Since Im $u_{n-1}$
is a torsion-free $R$-module, then $F_n$ is an RD-submodule of
$F_{n-1},$ which is not RD-injective. We deduce that $R$ is a
pure-injective module since it is a linearly compact module by
\cite[Proposition 9]{War1}, that fails to be RD-injective.}
\end{example}

\bigskip
It is proved in \cite[Theorem 3.1]{NaPr} that a domain $R$ is Pr\"ufer, i.e.
arithmetic, if and only if every pure-injective module is RD-injective.
This result can be extended to every commutative ring.
\begin{theorem} Let $R$ be a commutative ring. Then the following
assertions are equivalent~:
\begin{enumerate}
\item $R$ is an arithmetic ring.
\item Every RD-exact sequence of $R$-modules is pure-exact.
\item Every pure-injective $R$-module is RD-injective.
\end{enumerate}
\end{theorem}
\textbf{Proof.}
(1) $\Rightarrow$ (2) follows from \cite[Theorem 3]{War2}.
(2) $\Rightarrow$ (3) is obvious.
 (3) $\Rightarrow$ (1) is an immediate consequence of the following
proposition. \qed

\begin{proposition} Let $R$ be a commutative ring. If $R$ is not
arithmetic, there exists a pure-injective module which is not
RD-injective. More precisely, there exist a maximal ideal $P,$ and two
elements $a$ and $b$ in $P,$ such that $S = \{x\in E_R(R/P)\mid
ax=0 \hbox{ and } bx=0\}$ is a pure-injective module which is not RD-injective.
\end{proposition}
\textbf{Proof.}
There exists a maximal ideal $P$ such that $R_P$ is not a valuation ring.
Let $a$ and $b$ in $P$ such that $a\notin R_Pb$ and $b\notin R_Pa.$ By
\cite[Theorem 2]{War2}, there exists an indecomposable finitely presented
$R_P$-module
$M,$ generated by two elements $x_1$ and $x_2,$ with the relation $ax_1 =
bx_2.$ By \cite[Corollary 2]{War1}, $M$ is not an RD-projective $R_P$-module.
Consequently there exists an RD-exact sequence $0\rightarrow L\rightarrow
N\overset{\varphi}{\rightarrow} M\rightarrow 0$ which is not
pure-exact, where $N$ is an RD-projective $R_P$-module. We consider the
following presentation of $M:$  
 $0\rightarrow R_P\overset{u}{\rightarrow} F\overset
{v}{\rightarrow} M\rightarrow 0,$ where $F$ is a free $R_P$
-module with basis $\{e_1,e_2\},$ such that $v(e_i) = x_i,\ i=1,2$ and
$u(1) = ae_1-be_2.$ Since $F$ is a free module there exists $\alpha :
F\rightarrow N$ such that $\varphi\circ\alpha=v$ and consequently $\hbox{Im}(\alpha\circ u)\subseteq L.$
Then $a\alpha(e_1)-b\alpha(e_2)\notin aL+bL.$ Else there exist $y_1$
and $y_2$ in $L$ such that $ay_1-by_2 = a\alpha(e_1)-b\alpha(e_2).$ We 
define $\beta : F\rightarrow L$ by $\beta(e_i) = y_i.$ Then
$(\alpha-\beta)(ae_1-be_2) = 0,$ whence $(\alpha-\beta)$ induces
an homomorphism $\gamma : M\rightarrow N$ such that
$\varphi\circ\gamma =1_M.$
This is not possible.

Consequently if $G = R_P/(aR_P+bR_P),$ then
$G\otimes_{R_P} L \rightarrow G\otimes_{R_P}N$ is not injective. If
$E = E_R(R/P) \simeq E_{R_P}(R_P/PR_P),$ then
$E$ is an injective cogenerator in the category of $R_P$-modules, whence the
homomorphism
$\hbox{Hom}_{R_P}(G\otimes_{R_P}N,E)\rightarrow\hbox{Hom}_{R_P}(G\otimes_{R_P}L,E)$
is not surjective. If $S$ is the $R_P$-module $\hbox{Hom}_{R_P}(G,E),$
then $S$ is pure-injective by \cite[Proposition 7]{War1}, but the morphism
$\hbox{Hom}_{R_P}(N,S)\rightarrow\hbox{Hom}_{R_P}(L,S)$ is not
surjective, so $S$ is not RD-injective over $R_P.$

But $S \simeq
\hbox{Hom}_R(R/(aR+bR),E)\simeq \{x\in E\mid
ax=0\hbox{ and } bx=0\},$ and $S$ is also a pure-injective $R$-module that
fails to be RD-injective over $R.$ \qed

\bigskip
We deduce from this proposition the following example.

\begin{example} \textnormal{Let $K$ be a field and $R = K[X,Y]$ the polynomial ring
in two variables $X$ and $Y.$ If we take $a = X-\alpha,\ b = Y-\beta,\
\alpha,\beta\in K$ and $P$ the maximal ideal generated by $a$ and
$b,$ then $S$ is isomorphic to the simple $R$-module $R/P.$ This is
an example of a simple module which is not RD-injective over a domain
$R.$ Let us observe that $S$ is finite if $K$ is a finite field.}
\end{example}

\end{document}